\documentclass[12pt]{amsart}
\usepackage{euscript}
\usepackage{epsfig}
\usepackage{epic}
\usepackage{eepic}
\newcommand{\G}{\Gamma}
\newcommand{\PG}{P\G}
\newcommand{\E}{\mathcal{E}}
\renewcommand{\H}{\mathcal{H}}
\newcommand{\M}{\mathcal{M}}
\newcommand{\Z}{{\bf Z}}
\newcommand{\F}{{\bf F}}
\newcommand{\C}{{\bf C}}
\newcommand{\cp}{\C P}
\newcommand{\w}{\omega}
\newcommand{\wbar}{{\bar\omega}}
\newcommand{\thetabar}{{\bar\theta}}
\newcommand{\tensor}{\otimes}
\newcommand{\D}{\Delta}
\newcommand{\orbpi}{\pi^{\hbox{\rm\scriptsize orb}}}
\newcommand{\ip}[2]{\langle #1|#2\rangle}
\newcommand{\aut}{\mathop{\hbox{Aut}}\nolimits}
\newcommand{\Paut}{\mathop{P\kern-1pt\hbox{Aut}}\nolimits}
\newcommand{\diag}{\mathop{\hbox{diag}}\nolimits}
\newcommand{\im}{\mathop{\hbox{Im}}\nolimits}
\newcommand{\EE}[1]{\Lambda_4^\E} % so \EE8 gives Eisenstein E8
\renewcommand{\P}{{\bf P}}
\newcommand{\Eleech}{\Lambda_{12}^\E}
\newcommand{\GL}{{\rm GL}}
\newcommand{\suz}{\hbox{\it Suz}}
\newcommand{\sset}{\subseteq}
\newcommand{\isomorphism}{\cong}
\renewcommand{\a}{\alpha}
\newcommand{\DG}{\Lambda_2^{\mathcal{G}}} % gaussian D4
\renewcommand{\O}{{\rm O}} % orthogonal group

\theoremstyle{plain}
\newtheorem*{conjecture}{Conjecture}
\numberwithin{table}{section}
\numberwithin{figure}{section}

\begin{document}

\title{A Monstrous Proposal}
\author{Daniel Allcock}
\address{Department of Mathematics\\University of Texas at Austin\\Austin, TX 78712}
\email{allcock@math.utexas.edu}
\urladdr{http://www.math.utexas.edu/\textasciitilde allcock}
\date{June 30, 2007}
\thanks{Partly supported by NSF grant DMS-0245120.}
\dedicatory{Dedicated to Domingo Toledo and to John McKay}
\begin{abstract}
We explain a conjecture relating the monster simple group to an
algebraic variety that was discovered in a non-monstrous context.
\end{abstract}

\maketitle

\noindent
The purpose of this note is to explain a conjecture I have circulated
privately since 1997.  On one level it is purely about group theory
and complex hyperbolic geometry, but if it is true then the most
natural explanation for it would be algebra-geometric: a certain ball
quotient would be the moduli space for some sort of objects, the
objects admitting some sort of marking related to the monster simple
group.  My original grounds for the conjecture were flimsy, but in his
dissertation, Tathagata Basak  discovered some very suggestive
coincidences.  The conjecture is still speculative, but now I am
taking it seriously.

In brief, Conway described the bimonster $M\times M:2$ as being
generated by 16 involutions satisfying some braid and commutation
relations, subject to the additional relation that a certain word $w$
has order~10.  On the other hand, I discovered a certain complex
analytic orbifold $X$ having a nice uniformization by complex
hyperbolic 13-space $B^{13}$.  There is nothing obviously monstrous
about $X$, and I had no thought of the monster when constructing it.
Then I noticed that its fundamental group $\orbpi_1(X)$ contains 16
elements satisfying exactly the same braid and commutation relations,
but having order 3 rather than~2.  Later, Basak discovered that they
generate all of $\orbpi_1(X)$, that $w$ has order~20, and that the
generators can be extended to a set of 26 in exactly the same way as
Conway's bimonster generators.  This note formulates a conjecture that
explains these coincidences by giving a uniform interpretation of both
groups.

\section{Conjecture}
\label{sec-conjecture}

\noindent
In \cite{allcock}, we considered a certain lattice $L$ over the Eisenstein
integers $\E=\Z[\w{=}\sqrt[3]{1}]$, of signature $(1,13)$; a lattice
means a free module equipped with an $\E$-valued Hermitian form
$\ip{}{}$.
$L$ is the most natural such lattice, in the sense that its
underlying $\Z$-lattice is a scaled version of the unique even
unimodular lattice of signature $(2,26)$.  We will describe $L$
explicitly below, but for now it suffices to describe $L$ as the
unique $\E$-lattice of signature $(1,13)$ satisfying $L=\theta L^*$,
where $\theta=\w-\wbar=\sqrt{-3}$ and the asterisk denotes the dual
lattice.  In particular, all inner products in $L$ are divisible by
$\theta$.  Because of this, if $r\in L$ has norm $r^2=\ip{r}{r}=-3$
(such an $r$ is called a root of $L$), then the map
\begin{equation}
\label{eq-formula-for-reflection}
x\mapsto x-(1-\w)\frac{\ip{x}{r}}{\ip{r}{r}}r
\end{equation}
is an isometry of $L$; it sends $r$ to $\w r$ and fixes $r^\perp$
pointwise.   This is a complex
reflection of order~$3$, also called a triflection, and
$r^\perp$ is called its mirror.

The importance of $L$ in \cite{allcock} was that it allowed the
construction of a complex hyperbolic reflection group of record
dimension, namely the the subgroup $\G$ of $\aut L$ generated
by the triflections in the roots of $L$.  Now, $\aut L$  acts
on the complex $13$-ball $B^{13}$ consisting of positive lines in
$\C^{1,13}=L\tensor_\E\C$.  I proved that $\G$ has finite-volume
fundamental domain, and Basak has gone further \cite{basak-complex},
proving that
$\G$ is all of $\aut L$.

We define $X=B^{13}/\PG$, which is an algebraic variety and a complex
analytic orbifold.  We write $\H$ for the union of the mirrors of the
triflections, $\D$ for the image of $\H$ in $X$, and $G$
for the orbifold fundamental group $\orbpi_1(X-\D)$.  By a meridian we
mean an element of $G$ represented by a small loop in $X-\D$ that
encircles $\D$ once positively at a generic point of $\D$, or any
conjugate of such a loop.

\begin{conjecture}
The quotient of $G$ by the
normal subgroup generated by the squares of the meridians is
the bimonster, i.e., the semidirect product $M\times M:2$, where $M$
is the monster simple group and $\Z/2$ acts by exchanging the
factors.
\end{conjecture}

The conjecture gives a uniform description of the bimonster and $\PG$,
as the deck groups of the covering spaces of $X$ which are universal
among those having 2- and 3-fold branching along $\D$.  A geometric
way to understand it is to consider the orbifold structure along
$\D\sset X$.  Because of the triflections, a generic point of $\D$ has
local group $\Z/3$.  The idea of the conjecture is to rub out this
$\Z/3$ and replace it by $\Z/2$.  This should change the orbifold
fundamental group from $\PG$ to the bimonster.  The stated form of the
conjecture expresses this idea without having to make sense of the new
orbifold structure at non-generic points of $\D$.  (This can be done,
but involves a digression; see remark (6) in
section~\ref{sec-remarks}.)

If the conjecture holds then we should expect some reason for it to be
true.  One possibility that would be very pretty is for $X$ to have a
moduli interpretation.  If $\M$ is a moduli space parameterizing some
sort of algebra-geometric object, then a cover $\widetilde\M$ of $\M$,
possibly ramified somewhere, can be thought of as a moduli space
parameterizing suitably marked versions of the objects.

For example we may take $\M$ to be the moduli space of unordered
$12$-tuples in $\cp^1$ (say, in the sense of geometric invariant
theory) and $\widetilde\M$ to be the moduli space of ordered
$12$-tuples.  In this case the marking is the ordering, the deck group
is the symmetric group $S_{12}$, and $\widetilde\M$ is the universal
cover of $\M$ with $2$-fold branching along the discriminant.  This
example is especially relevant because $\M$ may also be described as a
$9$-ball quotient; indeed $B^9$ is the universal cover of
$\M-(\hbox{one point})$ with $3$-fold branching along the
discriminant.  So this situation is an exact analogue of that of the
conjecture, with $B^{13}$ replaced by $B^9$ and the bimonster by
$S_{12}$.  In fact, the group acting on $B^9$ is also the analogue of
our $\PG$.  See \cite{allcock}, \cite{DM} and \cite{thurston} for more
information.

It would be very pretty if the same phenomenon happened
for $X$; then the objects it parameterized would admit a sort of
marking, and varying the objects in a family would permute the
markings by an action of the bimonster.  The branched covering
$B^{13}\to X$ would parameterize the same objects, but equipped with a
different notion of marking.

\section{Evidence}
\label{sec-evidence}

\noindent
The origin of the conjecture is a coincidence of diagrams. In
\cite{conway-bimonster} Conway conjectured that the bimonster can be
presented as the quotient of the Coxeter group with diagram $Y_{555}$
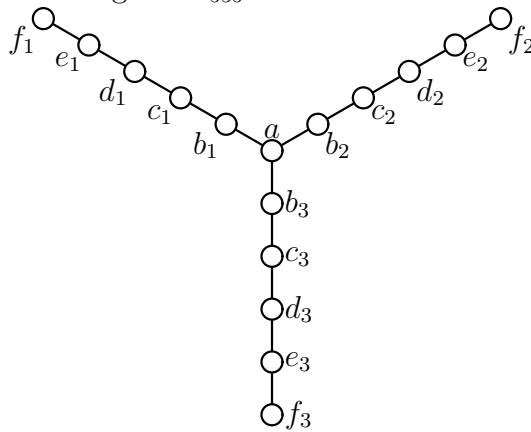
\begin{displaymath}
%
% the suffix D is a reminder that these variables should be dimens
\newdimen\edgelengthD  \edgelengthD=20pt
\newdimen\nodediameterD \nodediameterD=8pt
\newdimen\XshiftD \XshiftD=\edgelengthD \multiply\XshiftD by 173
\divide\XshiftD by 200
\newdimen\YshiftD \YshiftD=\edgelengthD \divide\YshiftD by 2
%
%
% U and V are points to use as local variables
\newcount\UxC
\newcount\UyC
\newcount\VxC
\newcount\VyC
\newcount\nodediameterC
% from (#1,#2) to (#3,#4) arguments should be dimensions
\def\bond#1#2#3#4{%
  \thicklines
  \UxC=#1
  \UyC=#2
  \VxC=#3
  \VyC=#4
  \drawline(\UxC,\UyC)(\VxC,\VyC)
}
% #1,#2 should be an ordered pair of dimensions
\def\hollownode#1#2{%
  \thicklines
  \UxC=#1
  \UyC=#2
  \nodediameterC=\nodediameterD
  \filltype{white}\put(\UxC,\UyC){\circle*{\nodediameterC}}
}
%
% FOR PLACING TEXT AND OTHER STUFF
%
% like ordinary put, except that args should be dimens
\def\myput#1#2#3{%
  \UxC=#1
  \UyC=#2
  \put(\UxC,\UyC){#3}
}
% (#1,#2) should be position (dimensions not numbers)
% (#3,#4) should be how many hundredths of a noderadius
% to offset the position.  (#5,#6) should be an additional
% offset (again, a dimension).  #7 should be the positioning
% string (tr, b, etc.)#8 should the the thing to typset.
%
\newcount\XoffC
\newcount\YoffC
\newdimen\XoffD
\newdimen\YoffD
\newdimen\noderadiusD
\newdimen\UxD
\newdimen\UyD
\def\nearnode#1#2#3#4#5#6#7#8{%
  \UxD=#1
  \UyD=#2
  \XoffC=#3
  \YoffC=#4
  \noderadiusD=\nodediameterD \divide\noderadiusD by 2
  \XoffD=\noderadiusD
  \YoffD=\noderadiusD
  \multiply\XoffD by \XoffC
  \multiply\YoffD by \YoffC
  \divide\XoffD by 100
  \divide\YoffD by 100
  \advance\UxD by \XoffD
  \advance\UyD by \YoffD
  \myput\UxD\UyD{\kern#5\makebox(0,0)[#7]{\raise#6\hbox{#8}}}
}
\setlength{\unitlength}{1sp}
\newcount\LLy
\newcount\height
\LLy=-\edgelengthD
\multiply\LLy by 5
\height=\edgelengthD
\multiply\height by 15
\divide\height by 2
\begin{picture}(0,\height)(0,\LLy)
%\put(0,\LLy){\framebox(0,\height)[bl]{\relax}}
% origin at A
%
%  Fone					Ftwo
%      Eone			    Etwo
%	   Done			Dtwo
%	       Cone	    Ctwo
%		   Bone	Btwo
%		       A
%		     Bthree
%
%		     Cthree
%
%		     Dthree
%
%		     Ethree
%
%		     Fthree
%
\newdimen\AX \AX=0pt
\newdimen\AY \AY=0pt
\newdimen\BoneX \BoneX=-\XshiftD
\newdimen\BoneY \BoneY=\YshiftD
\newdimen\ConeX \ConeX=-\XshiftD \multiply\ConeX by 2
\newdimen\ConeY \ConeY=\YshiftD \multiply\ConeY by 2
\newdimen\DoneX \DoneX=-\XshiftD \multiply\DoneX by 3
\newdimen\DoneY \DoneY=\YshiftD \multiply\DoneY by 3
\newdimen\EoneX \EoneX=-\XshiftD \multiply\EoneX by 4
\newdimen\EoneY \EoneY=\YshiftD \multiply\EoneY by 4
\newdimen\FoneX \FoneX=-\XshiftD \multiply\FoneX by 5
\newdimen\FoneY \FoneY=\YshiftD \multiply\FoneY by 5
\newdimen\BtwoX \BtwoX=\XshiftD
\newdimen\BtwoY \BtwoY=\YshiftD
\newdimen\CtwoX \CtwoX=\XshiftD \multiply\CtwoX by 2
\newdimen\CtwoY \CtwoY=\YshiftD \multiply\CtwoY by 2
\newdimen\DtwoX \DtwoX=\XshiftD \multiply\DtwoX by 3
\newdimen\DtwoY \DtwoY=\YshiftD \multiply\DtwoY by 3
\newdimen\EtwoX \EtwoX=\XshiftD \multiply\EtwoX by 4
\newdimen\EtwoY \EtwoY=\YshiftD \multiply\EtwoY by 4
\newdimen\FtwoX \FtwoX=\XshiftD \multiply\FtwoX by 5
\newdimen\FtwoY \FtwoY=\YshiftD \multiply\FtwoY by 5
\newdimen\BthreeX \BthreeX=0pt 
\newdimen\BthreeY \BthreeY=-\edgelengthD 
\newdimen\CthreeX \CthreeX=0pt
\newdimen\CthreeY \CthreeY=-\edgelengthD \multiply\CthreeY by 2
\newdimen\DthreeX \DthreeX=0pt
\newdimen\DthreeY \DthreeY=-\edgelengthD \multiply\DthreeY by 3
\newdimen\EthreeX \EthreeX=0pt
\newdimen\EthreeY \EthreeY=-\edgelengthD \multiply\EthreeY by 4
\newdimen\FthreeX \FthreeX=0pt
\newdimen\FthreeY \FthreeY=-\edgelengthD \multiply\FthreeY by 5
\bond\AX\AY\FoneX\FoneY
\bond\AX\AY\FtwoX\FtwoY
\bond\AX\AY\FthreeX\FthreeY
\hollownode\AX\AY
\hollownode\BoneX\BoneY
\hollownode\ConeX\ConeY
\hollownode\DoneX\DoneY
\hollownode\EoneX\EoneY
\hollownode\FoneX\FoneY
\hollownode\BtwoX\BtwoY
\hollownode\CtwoX\CtwoY
\hollownode\DtwoX\DtwoY
\hollownode\EtwoX\EtwoY
\hollownode\FtwoX\FtwoY
\hollownode\BthreeX\BthreeY
\hollownode\CthreeX\CthreeY
\hollownode\DthreeX\DthreeY
\hollownode\EthreeX\EthreeY
\hollownode\FthreeX\FthreeY
\nearnode\AX\AY{0}{120}{0pt}{0pt}{b}{$a$}%
\nearnode\BoneX\BoneY{-70}{-70}{0pt}{0pt}{tr}{$b_1$}%
\nearnode\ConeX\ConeY{-70}{-70}{0pt}{0pt}{tr}{$c_1$}%
\nearnode\DoneX\DoneY{-70}{-70}{0pt}{0pt}{tr}{$d_1$}%
\nearnode\EoneX\EoneY{-70}{-70}{0pt}{0pt}{tr}{$e_1$}%
\nearnode\FoneX\FoneY{-70}{-70}{0pt}{0pt}{tr}{$f_1$}%
\nearnode\BtwoX\BtwoY{70}{-70}{0pt}{0pt}{tl}{$b_2$}%
\nearnode\CtwoX\CtwoY{70}{-70}{0pt}{0pt}{tl}{$c_2$}%
\nearnode\DtwoX\DtwoY{70}{-70}{0pt}{0pt}{tl}{$d_2$}%
\nearnode\EtwoX\EtwoY{70}{-70}{0pt}{0pt}{tl}{$e_2$}%
\nearnode\FtwoX\FtwoY{70}{-70}{0pt}{0pt}{tl}{$f_2$}%
\nearnode\BthreeX\BthreeY{120}{0}{0pt}{0pt}{l}{$b_3$}%
\nearnode\CthreeX\CthreeY{120}{0}{0pt}{0pt}{l}{$c_3$}%
\nearnode\DthreeX\DthreeY{120}{0}{0pt}{0pt}{l}{$d_3$}%
\nearnode\EthreeX\EthreeY{120}{0}{0pt}{0pt}{l}{$e_3$}%
\nearnode\FthreeX\FthreeY{120}{0}{0pt}{0pt}{l}{$f_3$}%
\end{picture}
\end{displaymath}
by a single extra relation.  There are several ways to write the extra
relation, one of which is $(ab_1c_1ab_2c_2ab_3c_3)^{10}=1$, called the
spider relation.  Ivanov and Norton proved his conjecture in
\cite{ivanov} and \cite{norton}.

The diagram $Y_{550}$ arose as figure~5.1 of my paper \cite{allcock}, in a
manner suggesting its extension to $Y_{555}$.  It describes an
arrangement of~$11$ vectors in $\C^{1,9}$ of norm $-3$ such that their
triflections braid ($aba=bab$) or commute ($ab=ba$) according to
whether the corresponding nodes of the figure are joined or not.
The obvious generalization to $Y_{555}$ is the arrangement of
vectors in $\C^{1,13}$
$$
\renewcommand{\t}{$\theta$}
\renewcommand{\b}{$\thetabar$}
\newcommand{\1}{$1$}
\newcommand{\m}{$-1$}
\newcommand{\z}{$0$}
\newcommand{\W}{$\wbar$}
\newcommand{\0}{\phantom{$-1$}}
\begin{tabular}{l@{${}=($}r@{,}r@{,}r@{,}r@{;\,}r@{,}r@{,}r@{,}r@{;\,}r@{,}r@{,}r@{,}r@{;\,}r@{,}r@{)}}
$a  $&\0&\0&\0&\0&\0&\0&\0&\0&\0&\0&\0&\0&\1&\W\\
$b_1$&\0&\0&\0&\b&\0&\0&\0&\0&\0&\0&\0&\0&\z&\1\\
$c_1$&\0&\1&\1&\1&\0&\0&\0&\0&\0&\0&\0&\0&\0&\0\\
$d_1$&\0&\t&\0&\0&\0&\0&\0&\0&\0&\0&\0&\0&\0&\0\\
$e_1$&\m&\m&\1&\0&\0&\0&\0&\0&\0&\0&\0&\0&\0&\0\\
$f_1$&\b&\0&\0&\0&\0&\0&\0&\0&\0&\0&\0&\0&\0&\0\\
\end{tabular}
$$ with blanks indicating zeros and $b_2,\dots,f_2$ and
$b_3,\dots,f_3$ got from $b_1,\dots,f_1$ by permuting the first three
blocks of coordinates.  Here we are referring to the Hermitian form
with inner product matrix
$\diag[-1,\dots,-1]\oplus\bigl(\begin{smallmatrix}0&\thetabar\\\theta&0\end{smallmatrix}\bigr)$.
The triflections in these roots braid or commute according to the
$Y_{555}$ diagram, and the roots span a copy of $L$ in the form
$L=\EE8\oplus\EE8\oplus\EE8\oplus H$.  Here, $\EE8$ is the $E_8$ root
lattice regarded as a 4-dimensional $\E$-lattice and $H$ is the
``hyperbolic cell''
$\bigl(\begin{smallmatrix}0&\thetabar\\\theta&0\end{smallmatrix}\bigr)$.
For each $i=1,2,3$, the $i$th summand $\EE4$ is spanned by
$c_i,\dots,f_i$.

The Artin group $A$ of $Y_{555}$ is the group with one
generator for each node of the diagram, subject to the same
commutation and braiding relations as in the Coxeter group.  Forcing
the generators to have order~$2$ yields the Coxeter group of $Y_{555}$,
so Conway provides us with a surjection from $A$ to the bimonster.  On
the other hand, forcing them to have order~$3$ gives a map from $A$
into $\G$ by sending the Artin generators to the triflections in
$a$ and the $b_i,\dots,f_i$.  Basak has proven that these~$16$
triflections generate $\G=\aut L$, so both the bimonster and $\G$ are
quotients of $A$, with generators of order~$2$ and~$3$ respectively.
I expect that there is also a map $A\to G$, sending the Artin generators
to meridians, so that $A\to\PG$ is the composition
$A\to G\to\PG$.
%This would be very
%natural, because the quotient of a hyperplane complement by a group
%generated by reflections in the hyperplanes typically has
%Artin-group-like fundamental group.
%(I do not know whether $A\to G$ is surjective.)

This suggests that $G$ is really the central object, and that it
is `like' $A$, subject to some extra relations.  In the presence of
these extra relations, forcing the meridians to have order~$3$ reduces
the group to $\G$.  On the other hand, with luck, when the meridians
are forced to have order~$2$, the extra relations become equivalent to
the spider relation, giving the conjectured map $G\to M\times M:2$.  A
natural sanity check is whether the spider relation is compatible with
this picture.  Basak considered the word $S=ab_1c_1ab_2c_2ab_3c_3$ as
an element of $A$ and computed the order of its image in $\G$, which
turns out to be~$20$.  This implies that $S$, regarded as a element of
$G$, has order a multiple of~$20$ (or infinite order), which is
certainly compatible with the spider relation $S^{10}=1$ in the
bimonster.  Indeed, 20 is a notably small multiple of~10.  Basak also
found similar compatibilities using other words.

He has also made the striking discovery that the~$16$ triflections may
be extended to a set of~$26$ in exactly the same way as
Conway's~$16$ bimonster generators.
Namely, Conway observed that the map from the Coxeter group of
$Y_{555}$ to the bimonster extends to a map from the Coxeter group of
a larger graph, the incidence graph of the $13$ points and $13$ lines
of the projective plane over $\F_3$.  The symmetries of this finite
projective plane, including the dualities exchanging points and lines,
extend to automorphisms of the bimonster.  Basak found a set of $26$
roots of $L$, the reflections in which commute or braid according to
this graph; indeed all inner products are~$0$ except for
$\ip{p}{l}=\theta$ when $p$ is the root corresponding to a point of
$\P^2(\F_3)$ and $l$ is the root corresponding to a line containing
it.  It is not really surprising that the $16$ roots of $Y_{555}$ fit
in $\C^{1,13}$, but it {\it is\/} surprising that $26$ vectors with
specified inner products just happen to fit in a $14$-dimensional
space, and just happen to span a very natural lattice.  The mirrors in
$B^{13}$ corresponding to the $13$ points (resp. lines) of
$\P^2(\F_3)$ are mutually orthogonal and meet at a single point, say
$P$ (resp. $L$).  The midpoint of the segment joining $P$ and $L$ has
stabilizer $\GL_3(\F_3){:}2$ in $\PG$, realizing every automorphism of
$\P^2(\F_3)$ including dualities exchanging points and lines.
Presumably there is a map from the Artin group of the 26-node graph to
$G$, but I have not investigated this.

A final coincidence is that one of the maximal subgroups of the
monster has structure $3^{1+12}{\cdot}2{\cdot}\suz:2$, where $\suz$
denotes Suzuki's sporadic finite simple group, while $\PG$ contains a
subgroup $K$ with structure $(\im\E)\cdot\Eleech:(6{\cdot}\suz)$.  Here,
$\im\E\isomorphism\Z$ and $\Eleech$ denote the additive groups of the
imaginary part of $\E$ and of the complex Leech lattice.  The extension
defining the Heisenberg group $(\im\E)\cdot\Eleech$ is given by the
imaginary part of the inner product on $\Eleech$, and $6{\cdot}\suz$ is
$\aut\Eleech$.  Identifying the scalar $\w$ of $6{\cdot}\suz$ with a
generator of $\im\E$ reduces $K$ to $3^{1+12}{\cdot}2{\cdot}\suz$.
There are two ways to do this, corresponding to the two generators of
$\im\E$; this mimics the construction of
$3^{1+12}{\cdot}2{\cdot}\suz:2\sset M$ as a quotient of
$3^{1+12}:6{\cdot}\suz:2$ in \cite[secs. 3.3 and 5.2]{LPWW}.

$K$ is a very natural subgroup of $\PG$, namely the stabilizer of the
null vector $(0;0,1)$ in the realization of $L$ as
$\Eleech\oplus\bigl(\begin{smallmatrix}0&\thetabar\\\theta&0\end{smallmatrix}\bigr)$.
Since there are no mirrors in $B^{13}$ passing through $(0;0,1)$, $K$
is a subgroup of $G$, not just a subquotient of $G$.  Therefore we
have a natural subgroup of $G$, which modulo a natural relation has
index~2 in a maximal subgroup of the monster.

\section{Remarks}
\label{sec-remarks}

We close with some remarks that seem relevant but are not in the main
line of ideas.

(1)  There are infinitely many triflections in $\G$, but their mirrors
     form a locally finite arrangement in $B^{13}$ because $\G$ is
     discrete.  $\G$ acts transitively on roots, so all the
     triflections are conjugate (up to inversion)
     and $\D\sset X$ is irreducible.  The
     only complex reflections in $\G$ or $\PG$ are the triflections
     we have considered.

(2) There exist points of $B^{13}$ that lie on no mirrors, yet have
     nontrivial stabilizer in $\PG$.  For example, the stabilizer
     $6{\cdot}\suz$ is possible for points near the ideal point of
     $B^{13}$ given by
     $(0;0,1)\in\Eleech\oplus\bigl(
     \begin{smallmatrix}0&\thetabar\\\theta&0\end{smallmatrix}\bigr)$.
       A consequence is that $X-\D$ has nontrivial orbifold structure
    even though the obvious orbifold points have been removed.  This
    phenomenon does not occur for finite complex reflection groups or
    real hyperbolic Coxeter groups.

(3) Conway and Pritchard \cite{conway-fi24} used real
    hyperbolic geometry to study the bimonster and certain other
    finite groups as quotients of Coxeter groups.  I don't know of any
    connection between their real hyperbolic geometry and complex
    hyperbolic $13$-space  $B^{13}$.

(4) The study of complex reflection groups as quotients of Artin
    groups goes back to Coxeter \cite{coxeter}.  In our language, he
    studied the Artin groups for the diagrams $A_{n=1,\dots,5}$ modulo
    the relations that the generators have order~$3$.  In the first
    four cases one gets a finite complex reflection group acting on
    $\C^n$.  The last case gives a group of structure
    \hbox{$(\im\E)\cdot\EE8{:}\aut\EE8$}, which acts as a complex
    reflection group on the $5$-ball, fixing a point on the boundary.
    Coxeter showed that quotienting by the central $\im\E$ gives a
    complex reflection group acting cocompactly on $\C^4$.

(5) One can build up $Y_{555}$ from smaller diagrams by beginning with
    three $A_4$'s, which describe finite complex reflection groups,
    ``affinizing'' them by enlarging them to $A_5$'s (see the previous
    remark), and then ``hyperbolizing'' the result by adjoining a
    single extra node, joined to each of the three affinizing nodes.
    This is analogous to (say) extending the $E_8$ Coxeter group,
    which acts on the $7$-sphere, to $E_9$, which acts on Euclidean
    $8$-space, and then to $E_{10}$, which acts on real hyperbolic
    $9$-space.  Essentially the same process leads to diagrams like
    the one in \cite{borcherds}, which is a ``hyperbolization'' of
    $A_{11}D_7E_6$.

(6) We have avoided the issue of making sense of the alteration of
    orbifold structures described right after the conjecture.  But
    this dodge is not necessary.  Suppose $b\in B^{13}$, $S$ is its
    stabilizer in $\PG$, $R\sset S$ is the subgroup generated by the
    triflections which fix $b$, and $U$ is a small ball around $b$.
    Then $R$ is a direct product of some copies of the triflection
    groups associated to $A_1,\dots,A_4$ in remark (4).  These act on
    $U$ by the direct product of their triflection representations and
    possibly a subspace where all the factors act trivially.  By
    Chevalley's theorem, $U/R$ is a smooth variety.  Obviously, the
    image of $\H$ therein is a divisor $D$, and $S/R$ acts on $U/R$,
    preserving $D$.  Let $U'$ be the cover of $U/R$ which is universal
    among those having 2-fold ramification along $D$ and no other
    ramification, and let $R'$ be the deck group of $U'$ over $U/R$.
    By \cite{orlik-solomon}, $R'$ is the direct product of Coxeter
    groups associated to the same diagrams $A_1,\dots,A_4$ as for $R$,
    and $U'$ is an open set in $\C^{13}$.  Also, $R'$ acts on $U'$ as
    the product of the Coxeter groups' standard representations, again
    possibly with some fixed subspace.  Now, an element of $S/R$
    preserves $D$, so it has a lift to an automorphism of $U'$, in
    fact $|R'|$ many lifts.  We take $S'$ to be the group consisting
    of all such lifts.  Then
$$
U'/S'
\isomorphism (U'/R')\!\bigm/\!(S'/R')
\isomorphism (U/R)\!\bigm/\!(S/R)
\isomorphism U/S,
$$ where the first and third isomorphisms are of complex analytic
orbifolds.  The middle isomorphism is one of complex analytic
varieties, and is a complex analytic orbifold isomorphism away from
the image of $D\sset U/R=U'/R'$. We have equipped $U/S\sset
B^{13}/\PG$ with an orbifold structure in which the generic point of
$\D$ has local group $\Z/2$ rather than $\Z/3$.  It is easy to see
that such a structure is unique.

(7) The truth of the conjecture would imply that the orbifold of the
    previous remark is the quotient of a complex manifold by an action
    of the bimonster.  Then each component of the preimage of $\D$ is
    a smooth hypersurface fixed pointwise by an involution $\a$ in the
    bimonster.  The centralizer of $\a$ is $\langle\a\rangle$ itself
    times a copy of the monster, so the monster acts on the
    fixed-point set of $\a$.   We wonder if
    this complex $12$-manifold (or a suitable compactification of it)
    could serve as the $24$-dimensional ``monster manifold'' sought by
    Hirzebruch et.\ al.\ \cite[pp.~86--87]{hirzebruch} in their study
    of elliptic cohomology.

(8) We have suggested that $X$ may be a moduli space; we do know that
    it contains several lower-dimensional moduli spaces.  Namely, the
    moduli space of unordered $12$-tuples in $\cp^1$ is the quotient
    of $B^9$ by $\Paut\bigl(\EE8\oplus\EE8\oplus H\bigr)$, the moduli
    space of genus~$4$ curves is the quotient of $B^9$ by
    $\Paut\bigl(\EE8\oplus N\oplus [-3]\oplus H\bigr)$, where $N$ is
    the orthogonal complement of a root in $\EE8$, and the moduli
    space of cubic threefolds is the quotient of $B^{10}$ by
    $\Paut\bigl(\EE8\oplus\EE8\oplus[-3]\oplus H\bigr)$.  We remark
    that this last lattice may be constructed from $Y_{551}$ in the
    same manner as $L$ was from $Y_{555}$.  See \cite{DM},
    \cite{thurston}, \cite{allcock}, \cite{kondo}, \cite{act} and
    \cite{looijenga} for more information about these ball quotients.

(9) The Coxeter group of the $Y_{555}$ diagram appears in the
    monodromy of the $T_{666}$ surface singularity
    $x^6+y^6+z^6+\lambda xyz=0$, where $\lambda$ is a nonzero constant
    (see \cite[sec.~3.8]{AGV} and \cite{gabrielov}), and also in Mukai's analysis \cite{mukai} of the moduli
    space of $12$ points in $(\P^5)^5$.  I don't know know any way to
    fit these facts into my conjectural framework.  If there is a
    connection then Conway and Pritchard's real-hyperbolic constructions
    are probably also relevant; see remark~(3).

(10) Simons \cite{simons} has studied a simpler version of our situation,
    concerning $Y_{333}$ and a quotient of the Coxeter group with diagram the
    incidence graph of the $7$ lines and $7$ points of $\P^2(\F_2)$,
    which turns out to be $\O_8^-(2){:}2$.  I found $14$
    tetraflections (order $4$ complex reflections) of $B^7$ that
    satisfy the commutation and braid relations given by this diagram.
    One can choose the roots to have norm~$-2$ and inner products $0$
    or $1\pm i$, spanning the lattice $\DG\oplus\DG\oplus\DG\oplus
    \bigl(\begin{smallmatrix}0&1+i\\1-i&0\end{smallmatrix}\bigr)$,
    where $\mathcal{G}$ is the Gaussian integers $\Z[i]$ and $\DG$ is
    the $D_4$ root lattice regarded as a lattice over $\mathcal{G}$.
    Because these are tetraflections, it is easy to turn them into
    involutions: one simply reduces the lattice modulo $1+i$.  Then
    the tetraflections act on an $8$-dimensional vector space over
    $\F_2$ equipped with a quadratic form of minus type.  This provides a
    nice perspective on Simons' group, although it does not seem to give a new
    proof of his theorem.  Basak has developed these ideas in a quaternionic
    context in \cite{basak-quaternionic}.

\end{document}